\input amstex 
\documentstyle{amsppt}
\input bull-ppt
\keyedby{bull489/lic}
\define\Fix{\operatorname{Fix}}

\topmatter
\cvol{30}
\cvolyear{1994}
\cmonth{April}
\cyear{1994}
\cvolno{2}
\cpgs{212-214}
\ratitle
\title Dynamical Zeta Functions for Maps of the 
Interval\endtitle
\author David Ruelle\endauthor
\shortauthor{David Ruelle}
\shorttitle{Dynamical Zeta Functions for Maps of the 
Interval}
\address Institut des Hautes \'Etudes Scientifiques, 
Bures-sur-Yvette, France \endaddress
\date January 18, 1992\enddate
\subjclass Primary 58F20, 58F03; Secondary 
58F11\endsubjclass
\keywords Zeta function, transfer operator, topological 
pressure,
interval map\endkeywords
\abstract A dynamical zeta function $\zeta$ and a transfer 
operator
$\scr L$ are associated with a piecewise monotone map $f$ 
of the
interval $[0,1]$ and a weight function $g$. The analytic 
properties
of $\zeta$ and the spectral properties of $\scr L$ are 
related
by a theorem of Baladi and Keller under an assumption of 
``generating
partition''. It is shown here how to remove this 
assumption and, in
particular, extend the theorem of Baladi and Keller to the 
case when
$f$ has negative Schwarzian derivative.\endabstract
\endtopmatter

\document

Let $0=a_0<a_1<\cdots<a_N=1$. We write $X=[0,1]\subset\Bbb 
R$ 
and assume that $f$ is continuous $X\to X$ and strictly 
monotone
on the intervals $J_i=[a_{i-1},a_i]$. Furthermore, let 
$g\:X\to
\Bbb C$ have bounded variation. A {\it transfer operator\/}
$\scr L$ acting on functions $\Phi\: X\to \Bbb C$ of 
bounded variation
is defined by
$$
(\scr L\Phi)(x)=\sum_{y\:fy=x} g(y)\Phi(y),
$$
and we let
$$
\theta=\lim_{m\to\infty}\sup_{x\in X} \lf|\prod_{k=0}^{m-1}
g(f\,^kx)\rt|^{1/m}.
$$
It is known (see \cite1) that the essential spectral 
radius of
$\scr L$ is $\le\theta$. Assuming that $(J_1,\dots, J_N)$ is
generating (i.e., if $f\,^k x$ and $f\,^ky$ are in the same
$J_{i(k)}$ for all $k\ge 0$, then $x=y$\<), Baladi and 
Keller 
\cite1 have proved the following remarkable result 
(referred to
as B-K theorem in what follows):

\thm\nofrills{}\ The function
$$
d(z)=\exp -\sum_{m=1}^\infty\frac{z^m}m \sum_{x\in\Fix 
f\,^m}
\prod_{k=0}^{m-1} g(f\,^kx)
$$
is holomorphic for $|z|\theta<1$, and its zeros there are 
the
inverses $\lambda^{-1}$ of the eigenvalues $\lambda$ of 
$\scr L$
such that $|\lambda|>\theta$, with the same 
multiplicity.\ethm

Note that $\zeta(z)=1/d(z)$ is the {\it dynamical zeta 
function\/}
associated with the weight $g$ in the manner suggested by 
statistical
mechanics (see \cite8); Milnor and Thurston \cite7 have 
discussed the
case $g=1$. If $(J_1,\dots, J_N)$ is not generating, the 
set 
$\operatorname{Per} f$
of $f$-periodic points may be uncountable, so $\zeta(z)$ 
is not even
defined. The purpose of the present note is to indicate 
how to remove
this obstruction to the B-K theorem.

There is a map $\alpha$ defined on 
$\operatorname{Per}f\backslash\text{finite}$\ set such 
that $\alpha x=
(\xi(k))_{k\ge 0}$ is the unique sequence of symbols 
$\xi(k)\in\{1,\dots, N\}$
for which $f\,^kx\in J_{\xi(k)}$. We say that an 
$f$-invariant subset
$S$ of $\operatorname{Per}f$ is a {\it representative set 
of periodic
points\/} if the restriction of $\alpha$ to 
$S\backslash\text{finite}$\ 
set is a bijection to $(\alpha(\operatorname{Per}\backslash
\text{finite set}))\backslash\text{finite set}$\ and 
preserves the
period. Such a set $S$ always exists, and the B-K theorem 
remains true
provided $\zeta(z)$ is replaced by
$$
\zeta_S(z)=\exp \sum_{m=1}^\infty \frac{z^m}m
\sum_{x\in S\cap \Fix f\,^m} \prod_{k=0}^{m-1} g(f\,^kx).
$$
There are several cases of interest where 
$\operatorname{Per}f$ itself
is a representative set of periodic points: when 
$(J_1,\dots, J_N)$
is generating (case considered in \cite1), or when $f$ has 
negative
Schwarzian derivative: $Sf=f\,'''/f\,'-\tfrac 32 
(f\,''/f\,')^2<0$,
or when $f$ is affine on each $J_i$, and the slopes 
$\sigma_i$ are
such that $\prod\sigma_i^{m_i}\ne 1$ when $m_i\ge 0$, 
$\sum m_i>0$.
In those cases the B-K theorem remains true in its 
original form.

To prove these results, it is convenient to work in the 
more general 
setup where $X$ is a compact subset of $\Bbb R$ and to 
perform 
geometric changes on $(X,f,g,(J_1,\dots, J_N))$, observing 
the
effects of those changes on the spectrum of $\scr L$, 
$\theta$,
and on the zeta function. By doubling the points 
$a_1,\dots, 
a_{N-1}$ (see Hofbauer and Keller \cite5), collapsing 
certain
intervals (see Baladi and Ruelle \cite2), and embedding 
$(X,f\,)$
into a full shift on $N$ symbols (see Ruelle \cite 9), one 
reduces
to a particularly simple situation where one can use a 
form (due
to Baladi and Keller) of an argument due originally to 
Haydn \cite3.
This avoids the use of the Markov extension (an infinite 
Markov
partition) of Hofbauer \cite4 (see also \cite6). The 
effect of
the geometric changes mentioned above on $\theta$ and 
$\zeta$
is relatively easy to determine; the effect on the 
spectrum of
$\scr L$ takes more effort to analyze.

Using analytic completion (see Ruelle \cite9), one can 
extend the
analyticity properties of $d(z)$ to situations when $g$ is 
not a
function with bounded variation (but a quotient of such 
functions)
and one recovers results of Keller and Nowicki \cite6.

Let $r$ be the spectral radius of $\scr L$, and define the 
{\it pressure\/}
$$
P(\log|g|)=\sup_{\rho\in\bold I} \lf(h(\rho)+\int 
\rho(dx)\log |g(x)|\rt)
$$
where $h(\boldcdot)$ is the entropy of an element of the 
set $\bold I$
of $f$-invariant probability measures. With the above 
methods one can show
that
$$
\theta\le r\le\max(\theta,\exp P(\log|g|)).
$$
If $g\ge 0$ and $r>\theta$, then\fn"$^*$"{This result is 
obtained by
Baladi and Keller \cite1 with a different definition of 
the pressure.}
$r=\exp P(\log g)$ and
$r$ is an eigenvalue of $\scr L$; the set
$$
\Delta=\lf\{ \rho\in\bold I: h(\rho)+\int \rho(dx)\log 
g(x)=P(\log g)\rt\}
$$
of equilibrium measures is nonempty, $\Delta$ is a Choquet 
simplex,
and its vertices are ergodic measures.

\enddocument